\input amstex\documentstyle{amsppt}  
\pagewidth{12.5cm}\pageheight{19cm}\magnification\magstep1
\topmatter
\title Total positivity in symmetric spaces\endtitle
\author G. Lusztig\endauthor
\address{Department of Mathematics, M.I.T., Cambridge, MA 02139}\endaddress
\thanks{Supported by NSF grant DMS-1855773 and by a Simons Fellowship}\endthanks
\endtopmatter   
\document

\define\mpb{\medpagebreak}

\define\si{\sim}

\define\sqc{\sqcup}

\define\bK{\bar K}

\define\lb{\linebreak}

\define\part{\partial}

\define\m{\mapsto}
\define\do{\dots}

\define\lra{\leftrightarrow}

\define\sub{\subset}    

\define\T{\times}
\define\ti{\tilde}
\define\nl{\newline}
\redefine\i{^{-1}}
\define\fra{\frac}

\define\Hom{\text{\rm Hom}}

\define\a{\alpha}
\redefine\b{\beta}

\define\g{\gamma}
\redefine\d{\delta}
\define\e{\epsilon}

\redefine\o{\omega}

\define\ph{\phi}
\define\ps{\psi}
\define\r{\rho}
\define\s{\sigma}
\redefine\t{\tau}

\define\k{\kappa}

\define\z{\zeta}
\define\x{\xi}

\redefine\D{\Delta}

\define\Ps{\Psi}

\define\boc{\bold c}

\define\ii{\bold i}
\define\jj{\bold j}

\define\CC{\bold C}

\define\II{\bold I}

\define\NN{\bold N}

\define\QQ{\bold Q}
\define\RR{\bold R}

\define\ZZ{\bold Z}

\define\ci{\Cal I}
\define\cj{\Cal J}

\define\car{\Cal R}

\define\cu{\Cal U}

\define\fg{\frak g}
\define\fh{\frak h}

\define\ft{\frak t}

\define\tu{\ti u}

\define\bul{\bullet}

\subhead 0.1\endsubhead
Let $K$ be a real closed field that is, a field of characteristic zero
such that $K$ is not algebraically closed but $\bK:=K[\sqrt{-1}]$ is
algebraically closed.
Let $G$ be a connected reductive algebraic
group over $\bK$ with a pinning
$(B^+,B^-,T,x_i,y_i (i\in I))$ as in \cite{L94, 1.1}.
Let $U^+,U^-$ be the unipotent radicals of $B^+,B^-$.
Recall that $x_i:\bK@>>>U^+,y_i:\bK@>>>U^-$ are certain imbeddings
of algebraic groups. We identify $G$ with its group of $\bK$-points.

Let $\s:G@>>>G$ be the antiautomorphism of $G$ such that
$$\s(x_i(a))=x_i(a),\s(y_i(a))=y_i(a)$$
for all $i\in I,a\in\bK$ and $\s(t)=t\i$ for all $t\in T$.
Let $\o:G@>>>G$ be an involutive automorphism of $G$ preserving $T$ and
such that for some involution $i\m i^*$ of $I$ we have
$$\o(x_i(a))=x_{i^*}(a), \o(y_i(a))=y_{i^*}(a)$$ for all $i\in I,a\in\bK$.
We have $\s^2=1$, $\s\o=\o\s$. We set $\t=\s\o=\o\s$; this is an involutive
antiautomorphism of $G$ preserving $T$ and such that
$$\t(x_i(a))=x_{i^*}(a), \t(y_i(a))=y_{i^*}(a)$$ for all $i\in I,a\in\bK$.
Hence $g\m \t(g\i)$ is an involutive automorphism of $G$. Let
$H=\{g\in G;\t(g\i)=g\}$ be the fixed point set of this automorphism; then
$G/H$ is a symmetric space (of a special type).
Let $G^\t=\{g\in G;\t(g)=g\}$ (a closed smooth subvariety of $G$).

Now $G$ acts on $G$ by

(a) $g:g_1\m gg_1\t(g)$.
\nl
The stabilizer of $1$ in this action is $H$ so that $G/H$ can be
identified with the
$G$-orbit $(G^\t)^0$ of $1$ under this action via $gH\m g\t(g)$.
We have $(G^\t)^0\sub G^\t$. (This containment can be strict.)

\subhead 0.2\endsubhead
Since $U^+,U^-$ are $\t$-stable,
$U^{+\t}=U^+\cap G^\t$ is a closed smooth subvariety of $U^+$
on which $U^+$ acts by 0.1(a) and $U^{-\t}=U^-\cap G^\t$ is a closed
smooth subvariety of $U^-$ on which $U^-$ acts by 0.1(a).

\subhead 0.3\endsubhead
Let $G(K)$ be the group of $K$-points of $G$, a subgroup of $G$.
Let $G^\t(K)=G^\t\cap G(K)$,
$(G^\t)^0(K)=(G^\t)^0\cap G(K)$, $U^{+\t}(K)=U^{+\t}\cap G(K)$,
$U^{-\t}(K)=U^{-\t}\cap G(K)$,
$T(K)=T\cap G(K)$. Let $(G/H)(K)$ be the subset of $G/H$ corresponding to 
$(G^\t)^0(K)$ under the identification $(G^\t)^0=G/H$.

\subhead 0.4\endsubhead
Let $K_{>0}=\{a^2;a\in K-\{0\}\}\sub K-\{0\}$.
Then $K_{>0}$ is closed under addition, multiplication and division
(it is a semifield). Let $K_{\ge0}=K_{>0}\cup\{0\}\sub K$.
For $a,b$ in $K$ we write $a>b$ (resp. $a\ge b$) whenever $a-b\in K_{>0}$
(resp. $a-b\in K_{\ge0}$); this is a total order on $K$.

Let $T_{>0}=\{t^2;t\in T(K)\}$, a subgroup of $T(K)$. Following \cite{L94} let

$U^+_{\ge0}$ be the submonoid of $U^+$ generated by $x_i(a)$ with
$i\in I,a\in K_{\ge0}$;

$U^-_{\ge0}$ be the submonoid of $U^-$ generated by $y_i(a)$ with
$i\in I,a\in K_{\ge0}$;

$G_{\ge0}$ be the submonoid of $G$ generated by $x_i(a),y_i(a)$ with
$i\in I,a\in K_{\ge0}$ and by $T_{>0}$.

We set $G^\t_{\ge0}=G^\t\cap G_{\ge0}$, $U^{+\t}_{\ge0}=U^{+\t}\cap U^+_{\ge0}$,
$U^{-\t}_{\ge0}=U^{-\t}\cap U^-_{\ge0}$.

The subsets $G^\t_{\ge0},U^{+\t}_{\ge0},U^{-\t}_{\ge0}$ are not monoids but instead are sets
with action of the monoid $G_{\ge0},U^+_{\ge0},U^-_{\ge0}$ (respectively) which act by 0.1(a).

In 3.7 we show that $G^\t_{\ge0}$ is actually a subset of $(G^\t)^0(K)$
hence it can be viewed
as a subset of $(G/H)(K)$ said to be the totally positive part of $(G/H)(K)$.

In this paper we show that several properties of the
totally positive monoids $G_{\ge0},U^+_{\ge0},U^-_{\ge0}$ discussed in
\cite{L94}, such as existence of natural ``cell'' decompositions  (1.2(d), 1.2(j))
have analogues for the totally positive subsets
$G^\t_{\ge0},U^{+\t}_{\ge0},U^{-\t}_{\ge0}$ of $G^\t,U^{+\t},U^{-\t}$.
(See 2.4(b), 2.5(c), 3.8(b). We will write ``cell'' for something that is really a
cell when $K=\RR$.)
Actually most of the results in this paper are not only analogous to results in \cite{L94}
but contain those results in \cite{L94} as a special case (see 5.3).
But there is one property of cells in \cite{L94} which needs to be modified in the present case.

Let $W$ be the Weyl group of $G$. In \cite{L94} to each $w\in W$ we associate a cell
parametrized in several ways by $K_{>0}^n$ (with $K=\RR$) in such a way that
one obtains a positive structure on that cell, that is any two parametrizations
are related by a transformation given by substraction free rational functions.
One of our results is that for any twisted involution (see 2.1) $w\in W$ one can
attach a cell parametrized in several ways by $K_{>0}^{n'}$ (with $K=\RR$) in such a way that
(with some restriction on $*$) one obtains something more general than a positive
structure in the sense that
any two parametrizations are related by a transformation given by substraction free
rational functions
combined with extraction of square roots. See \S6.
We regard this as the main contribution of this paper.

\subhead 0.5\endsubhead
Many results of this paper extend (with essentially the same proof) to the case where
$G$ is replaced by a group attached to a crystallographic, possibly not positive definite,
symmetric Cartan matrix.

{\it Notation.} For any $c\in K_{>0}$ we denote by $\sqrt{c}$ the $>0$ square root of $c$.

\subhead 0.6\endsubhead
I thank Xuhua He for useful discussions.

\head 1. Preliminaries on $G_{\ge0},U^+_{\ge0},U^-_{\ge0}$\endhead
\subhead 1.1\endsubhead
For $i\in I$ let $s_i\in W$ be the simple reflection
corresponding to $i$. Let $S=\{s_i;i\in I\}$. Let $w\m|w|$ be the length function on $W$.
There is a well defined monoid structure $w,w'\m w\bul w'$ on $W$
characterized by

$s_i\bul w=s_iw$ if $i\in I,w\in W, |s_iw|=|w|+1$,

$s_i\bul w=w$ if $i\in I,w\in W, |s_iw|=|w|-1$.
\nl
For $w\in W$ let $\ci(w)$ be the set of sequences
$(i_1,i_2\do,i_k)$ in $I$ such that $k=|w|$ and
$w=s_{i_1}s_{i_2}\do s_{i_k}$.

Let $w_0$ be the longest element of $W$.

\subhead 1.2\endsubhead
We now recall some definitions and results from \cite{L94}.
Let $w\in W$. For $\ii=(i_1,i_2,\do,i_k)\in\ci(w)$ we define 
$\Ps^+_\ii:(K-\{0\})^k@>>>U^+$ by

$(a_1,a_2,\do,a_k)\m x_{i_1}(a_1)x_{i_2}(a_2)\do x_{i_k}(a_k)$
\nl
and $\Ps^-_\ii:(K-\{0\})^k@>>>U^-$ by

$(a_1,a_2,\do,a_k)\m y_{i_1}(a_1)y_{i_2}(a_2)\do y_{i_k}(a_k)$.
\nl
By the proof of \cite{L94, 2.7(a)}, 

(a) $\Ps^+_\ii,\Ps^-_\ii$ are injective.
\nl
Let $\ps^+_\ii:K_{>0}^k@>>>U^+$, $\ps^-_\ii:K_{>0}^k@>>>U^-$ be the
restrictions of $\Ps^+_\ii,\Ps^-_\ii$. By \cite{L94},

(b) $\ps^+_\ii$ defines a bijection of $K_{>0}^k$ onto
a subspace $U^+_{\ge0,w}$ of $U^+_{\ge0}$ independent of the
choice of $\ii$ in $\ci(w)$;

(c) $\ps^-_\ii$ defines a bijection of $K_{>0}^k$ onto
a subspace $U^-_{\ge0,w}$ of $U^-_{\ge0}$ independent of the
choice of $\ii$ in $\ci(w)$;

(d) we have a partition

$U^+_{\ge0}=\sqc_{w\in W}U^+_{\ge0,w}$;
\nl
we have a partition

$U^-_{\ge0}=\sqc_{w\in W}U^-_{\ge0,w}$.
\nl
The subsets $U^+_{\ge0,w}$ of $U^+_{\ge0}$ are said to be the ``cells'' of $U^+_{\ge0}$;
the subsets $U^-_{\ge0,w}$ of $U^-_{\ge0}$ are said to be the ``cells'' of $U^-_{\ge0}$.

For any $w,w'$ in $W$ we set

(e) $G_{\ge0,w,w'}=U^-_{\ge0,w}T_{>0}U^+_{\ge0,w'}=U^+_{\ge0,w'}T_{>0}U^-_{\ge0,w}
\sub G_{\ge0}$.
\nl
By \cite{L94}, we have bijections

(f) $U^-_{\ge0,w}\T T_{>0}\T U^+_{\ge0,w'}@>\si>>G_{\ge0,w,w'}$,

(g) $U^+_{\ge0,w'}\T T_{>0}\T U^-_{\ge0,w}@>\si>>G_{\ge0,w,w'}$,
\nl
given by multiplication in $G$. By \cite{L94},

(i) we have a partition $G_{\ge0}=\sqc_{(w,w')\in W\T W}G_{\ge0,w,w'}$.
\nl
The subsets $G_{\ge0,w,w'}$ of $G_{\ge0}$ are sad to be the ``cells'' of $G_{\ge0}$.

\subhead 1.3\endsubhead
By \cite{L94}, for $w,w'$ in $W$, multiplication in $G$ satisfies

(a) $U^+_{\ge0,w}U^+_{\ge0,w'}\sub U^+_{\ge0,w\bul w'}$,

(b) $U^-_{\ge0,w}U^-_{\ge0,w'}\sub U^-_{\ge0,w\bul w'}$.
\nl
It follows that for $w_1,w_2,w'_1,w'_2$ in $W$, multiplication
in $G$ satisfies

(c) $G_{\ge0,w_1,w'_1}G_{\ge0,w_2,w'_2}\sub
G_{\ge0,w_1\bul w_2,w'_1\bul w'_2}$.

\subhead 1.4\endsubhead
In this subsection we fix $w\in W,i\in I$ such that
setting $s_i=s$, we have $s^*w=ws,|w|=|s^*w|+1$. Let
$\ii'=(i_1,i_2,\do,i_k)\in\ci(s^*w)$,
$\ii=(i^*,i_1,i_2,\do,i_k)\in\ci(w),\ti\ii=(i_1,i_2,\do,i_k,i)\in\ci(w)$.
We define

(a) $\a:U^+_{\ge0,s^*w}\T K_{>0}@>>>U^+_{\ge0,w}$
\nl
by $(u,a)\m x_{i^*}(a)ux_i(a)$. This is well defined since
$s^*\bul(s^*w)\bul s=w$ (see 1.3(a)). We show:

(b) {\it The map $\a$ is a bijection.}
\nl
We first show that $\a$ is injective. Assume that $u,u'$
in $U^+_{\ge0,s^*w}$ and $a,a'$ in $K_{>0}$ satisfy
$x_{i^*}(a)ux_i(a)=x_{i^*}(a')u'x_i(a')$. Then $x_{i^*}(a-a')u=u'x_i(a'-a)$.

If $a-a'>0$, then by 1.2(b) we have

$x_{i^*}(a-a')u\in \ps^+_\ii(K_{>0}^{k+1})=\ps^+_{\ti\ii}(K_{>0}^{k+1})$
\nl
hence $x_{i^*}(a-a')u=u_1x_i(b)$ where $u_1\in U^+_{\ge0,s^*w}$ and
$b\in K_{>0}$. Thus we have $u_1x_i(b)=u'x_i(a'-a)$.
Using 1.2(a) (recall that $a'-a\ne0$) we see that
$b=a'-a$. Thus, $a'-a>0$, a contradiction.

If $a'-a>0$, then by 1.2(b) we have

$u'x_i(a'-a)\in \ps^+_{\ti\ii}(K_{>0}^{k+1})=
\ps^+_\ii(K_{>0}^{k+1})$
\nl
hence $u'x_i(a'-a)=x_{i^*}(b)u_1$
where $u_1\in U^+_{\ge0,s^*w}$ and $b\in K_{>0}$.
Thus we have $x_{i^*}(b)u_1=x_{i^*}(a-a')u$.
Using 1.2(a) (recall that $a-a'\ne0$) we see that
$b=a-a'$. Thus $a-a'>0$, a contradiction.

We see that we must have $a=a'$. It follows that $u=u'$
and the injectivity of $\a$ is proved.

We now prove that $\a$ is surjective. 
It is enough to show that for any $u\in U^+_{\ge0,s^*w}$ and
any $b\in K_{>0}$ we have $x_{i^*}(b)u=x_{i^*}(b')u'x_i(b')$
for some $u'\in U^+_{\ge0,s^*w}$ and some $b'\in K_{>0}$.
This is proved by the following argument, inspired by one
which I have learned from Xuhua He. \footnote{ X. He showed, in connection
with a different problem, that 
$\a$ is surjective when $K=\RR$, $W$ is of type $B_2$ and $|w|=4$.}

We can write uniquely
$u=x_{i_1}(a_1)x_{i_2}(a_2)\do x_{i_k}(a_k)$
with $a_1,a_2\do,a_k$ in $K_{>0}$.
For $c\in K_{>0}$ we have
$x_{i^*}(c)u=x_{i^*}(c)x_{i_1}(a_1)x_{i_2}(a_2)\do x_{i_k}(a_k)$
Using 1.2(b) we see that we have also
$x_{i^*}(c)u=x_{i_1}(b_1)x_{i_2}(b_2)\do x_{i_k}(b_k)x_i(b_{k+1})$
where $b_1,b_2,\do,b_{k+1}$ are uniquely determined in $K_{>0}$.
We shall regard $b_1,b_2,\do,b_{k+1}$ as functions of $c$ (here $a_1,a_2,\do,a_k$
are fixed). By results in \cite{L94}, (see also \cite{L19}), each of
$b_1,b_2,\do,b_{k+1}$
can be expressed as a substraction free rational function in $c,a_1,a_2,\do,a_k$.
In particular, there exist nonzero polynomials $P(\e),Q(\e)$ in $K[\e]$ ($\e$ is an indeterminate)
with all coefficients in $K_{\ge0}$ such that $b_{k+1}=P(c)/Q(c)$ when $c\in K_{>0}$; note that
$Q(c)>0$ when $c\in K_{>0}$. (Similar expressions
exist for the other $b_i$ but we do not need them). We can assume that not both $P(\e),Q(\e)$
are divisible by $\e$.

Now one checks easily that the map $(\bK-\{0\})^k\T\bK@>>>U^+$
given by $$(d_1,d_2,\do,d_k,d_{k+1})\m x_{i_1}(d_1)x_{i_2}(d_2)\do x_{i_k}(d_k)x_i(d_{k+1})$$
is an isomorphism of $(\bK-\{0\})^k\T\bK$ onto a subvariety $Z$ of $U^+$.
Let $\g:Z@>>>(\bK-\{0\})^k\T\bK$ be the inverse of this isomorphism.
Note that $c\m x_{i^*}(c)u$ is a morphism $\d:\bK@>>>U^+$ such that $\d(K_{\ge0})\sub Z$.
(We have $\d(0)=u\in Z$.) It follows that $\d\i(Z)$ is a subvariety of $\bK$ containing
$K_{\ge0}$ hence there exists an open subset $\cu$ of $\bK$ such that $K_{\ge0}\sub\cu$
and $\d$ restricts to a morphism $\d':\cu@>>>Z$.
Then the composition
of $\g\d':\cu@>>>(\bK-\{0\})^k\T\bK$ with the last projection $(\bK-\{0\})^k\T\bK@>>>\bK$
is a morphism $\r:\cu@>>>\bK$. In other words $c\m b_{k+1}$ from $K_{>0}$ to $K_{>0}$ is the
restriction of the morphism $\r:\cu@>>>\bK$. We have $\r(0)=0$.
(Indeed, we have  $\d(0)=u=\g(a_1,\do,a_k,0)$.) We see that $c\m P(c)/Q(c)$ from $K_{>0}$ to
$K_{>0}$ is the restriction of $\r:\cu@>>>\bK$ and $\r(0)=0$.
It follows that $P(0)=0,Q(0)\ne0$ and that $\r(c)=P(c)/Q(c)$ for $c\in K_{\ge0}$.
Since all coefficients of $Q(\e)$ are in $K_{\ge0}$ it follows that $Q(0)>0$.

For $c\in K_{\ge0}$ we have
$x_{i^*}(c)u=u'(c)x_i(P(c)/Q(c))$ where $u'(c)\in U^+_{\ge0,s^*w}$.

Let $[0,b]=\{b'\in K;0\le b'\le b\}$. If $b'\in[0,b]$ we have
$$\align&x_{i^*}(b)u=x_{i^*}(b')x_{i^*}(b-b')u=x_{i^*}(b')u'(b-b')x_i(P(b-b')/Q(b-b'))\\&
=x_{i^*}(b')u'(b-b')x_i(\fra{P(b-b')}{Q(b-b')}-b')x_i(b').\endalign$$
Next we note that 

(c) the function $b'\m r(b'):=P(b-b')-b'Q(b-b')$ from $[0,b]$ to $K$
takes the value $0$ for some $b'\in K,0<b'<b$.
\nl
This follows from the intermediate value theorem (known to hold for our $K$) applied
to the polynomial function $r:[0,b]@>>>K$. This function changes sign: we have
$r(0)=P(b)>0$ and $r(b)=P(0)-bQ(0)=-bQ(0)<0$.

For $b'$ as in (c) we have
$$x_{i^*}(b)u=x_{i^*}(b')u'(b-b')x_i(b').$$
This completes the proof of (b).

\subhead 1.5\endsubhead
In the setup of 1.4 recall that
$$\ii'=(i_1,i_2,\do,i_k)\in\ci(s^*w),\ti\ii=(i_1,i_2,\do,i_k,i)\in\ci(w).$$
We identify $U^+_{\ge0,s^*w}=K_{>0}^k$ via $\ps_{\ii'}$ and
$U^+_{\ge0,w}=K_{>0})^{k+1}$ via $\ps_{\ti\ii}$. Then $\a$ in 1.4 becomes
a bijection $K_{>0}^k\T K_{>0}@>>>K_{>0}^{k+1}$.
One can show that the inverse of this bijection is not in general
given by rational
functions. For example, if $W$ is of type $A_3$, $w=w_0$, $*=1$ and $s\in S$ is such that
$sw=ws$ then the formula for this inverse involves taking roots of quadratic
equations.

\subhead 1.6\endsubhead
In this subsection we fix $w\in W, w'\in W,i\in I$ such that,
setting $s_i=s$, we have $s^*w=ws,|w|=|s^*w|+1$. We define

(a) $\a:G_{\ge0,s^*w,w'}\T K_{>0}@>>>G_{\ge0,w,w'}$
\nl
by $(g,a)\m y_{i^*}(a)gy_i(a)$. This is well defined since
$s^*\bul(s^*w)\bul s=w$ (see 1.3(c)). We show:

(b) {\it the map $\a$ is a bijection. }
\nl
The proof is almost a copy of that of 1.4.
We first show that  $\a$ is injective. Assume that
$g,g'$ in $G_{\ge0,s^*w,w'}$ and $a,a'$ in $K_{>0}$ satisfy
$y_{i^*}(a)gy_i(a)=y_{i^*}(a')g'y_i(a')$. Then $y_{i^*}(a-a')g=g'y_i(a'-a)$.

If $a-a'>0$ then using 1.2(c),(e) we see that
$y_{i^*}(a-a')g=g_1y_i(b)$ where $g_1\in G_{\ge0,s^*w,w'}$ and
$b\in K_{>0}$. Thus we have $g_1y_i(b)=g'y_i(a'-a)$.
We can write $g_1=u_1t\tu_1$, $g'=u't'\tu'$
with $u_1,u'$ in $U^+_{\ge0,w'}$, $\tu_1,\tu'$ in
$U^-_{\ge0,s^*w}$, $t,t'$ in $T_{>0}$. From 
$u_1t\tu_1y_i(b)=u't'\tu'y_i(a'-a)$ we deduce
$u_1=u',t=t'$, $\tu_1y_i(b)=\tu'y_i(a'-a)$ (uniqueness of Bruhat
decomposition). From the last equality we see using 1.2(b) and
$a'-a\ne0$ that $b=a'-a$. Thus $a'-a>0$, a contradiction.

If $a'-a>0$ then using 1.2(c),(e) we see that
$g'y_i(a'-a)=y_{i^*}(b)g_1$
where $g_1\in G_{\ge0,s^*w,w'}$ and
$b\in K_{>0}$. Thus we have $y_{i^*}(b)g_1=y_{i^*}(a-a')g$.
We can write $g_1=\tu_1t u_1$, $g=\tu't'u'$
with $u_1,u'$ in $U^+_{\ge0,w'}$, $\tu_1,\tu'$ in
$U^-_{\ge0,s^*w}$, $t,t'$ in $T_{>0}$. From 
$y_{i^*}(b)u_1t\tu_1=y_{i^*}(a-a')u't'\tu'$ we deduce
$u_1=u',t=t'$, $y_{i^*}(b)u_1=y_{i^*}(a-a')u'$ (uniqueness of Bruhat
decomposition). From the last equality we see using 1.2(b) and
$a-a'\ne0$ that $b=a-a'$. Thus $a-a'>0$, a contradiction.

We see that we must have $a=a'$. It follows that $g=g'$
and the injectivity of $\a$ is proved.

We now prove that $\a$ is surjective.
It is enough to show that for any $g\in G_{\ge0,sw,w'}$ and
any $b\in K_{>0}$ we have $y_{i^*}(b)g=y_{i^*}(b')g'y_i(b')$
for some $g'\in G_{\ge0,s^*w}$ and some $b'\in K_{>0}$.

Using 1.2(e) we see that there  are well defined maps
$c\to z(c),K_{\ge0}@>>>K_{\ge0}$ and $c\to g'(c)$,
$K_{\ge0}@>>>G_{\ge0,s^*w,w'}$ 
such that $y_{i^*}(c)g=g'(c)y_i(z(c))$ for all $c\in K_{>0}$. We have $z(0)=0$, $g'(0)=g$.
As in 1.4 we see that there exist nonzero polynomials $P(\e),Q(\e)$ in $K[\e]$ with all
coefficients in $K_{\ge0}$ such that $Q(c)\ne0,z(c)=P(c)/Q(c)$ when $c\in K_{>0}$.
As in 1.4 we see that $c\m z(c)$ is the restriction to $K_{>0}$ of a morphism
$\r$ from an open subset of $\bK$ containing $K_{\ge0}$ to $\bK$ and that $\r(0)=0$.
As in 1.4 we deduce that we can assume $P(0)=0,Q(0)>0$. 
Let $[0,b]$ be as in 1.4. If $b'\in[0,b]$ we have
$$\align&y_{i^*}(b)g=y_{i^*}(b')y_i(b-b')g=y_{i^*}(b')g'(b-b')y_i(z(b-b'))\\&
y_{i^*}(b')g'(b-b')y_i(\fra{P(b-b')}{Q(b-b')}-b')y_i(b').\endalign$$
As in 1.4 we see using the intermediate value theorem for $K$ that for some
$b'\in K,0<b'<b$ we have 
$\fra{P(b-b')}{Q(b-b')}-b'=0$
hence
$$y_{i^*}(b)g=y_{i^*}(b')g'(b-b')y_i(b').$$
This completes the proof of (b).

\subhead 1.7\endsubhead
In this subsection we fix $w,w'$ in $W$, $i\in I$ such that,
setting $s_i=s$, we have $s^*w'=w's,|w'|=|s^*w'|+1$. We define

(a) $\a:G_{\ge0,w,s^*w'}\T K_{>0}@>>>G_{\ge0,w,w'}$
\nl
by $(g,a)\m x_{i^*}(a)gx_i(a)$. This is well defined since
$s^*\bul(s^*w')\bul s=w'$ (see 1.3(c)). The proof of the
following result is entirely similar to that of 1.6(b):

(b) {\it the map $\a$ is a bijection.}     

\subhead 1.8\endsubhead
Let $w,w'$ in $W$ and $i\in I$ be such that setting
$s=s_i\in S$ we have $|s^*ws|=|w|-2$. We show:

(a) {\it If $g,g'$ in $G_{\ge0,s^*ws,w'}$ and $a,a',c,c'$ in $K_{>0}$
satisfy  $$y_{i^*}(a)gy_i(c)=y_{i^*}(a')g'y_i(c')$$ then $a=a'$.}
\nl
We have $g=ut\tu$, $g'=u't'\tu'$ where
$u,u'$ are in $U^-_{\ge0,s^*ws}$, $\tu,\tu'$ are in $U^+_{\ge0,w'}$
and $t,t'$ are in $T_{>0}$. We have $$y_{i^*}(a)gy_i(c)=y_{i^*}(a)uy_i(b)t_1\tu_1,$$
$$y_{i^*}(a')g'y_i(c')=y_{i^*}(a')u'y_i(b')t'_1\tu'_1,$$
where $\tu_1,\tu'_1$ are in $U^+_{\ge0,w'}$,
$t_1,t'_1$ are in $T_{>0}$ and $b,b'$ are in $K_{>0}$.
We have
$$y_{i^*}(a)uy_i(b)t_1\tu_1=y_{i^*}(a')u'y_i(b')t'_1\tu'_1.$$
By the uniqueness in Bruhat decomposition we deduce
$$y_{i^*}(a)uy_i(b)=y_{i^*}(a')u'y_i(b').$$
Using this together with 1.2(a) and with $|w|=|s^*ws|+2$
we deduce that $a=a'$ proving (a).

\subhead 1.9\endsubhead
Let $w\in W$ with $|w|=l$. For any $\ii,\ii'$ in $\ci(w)$ we set
$u_{\ii,\ii'}=\ps_{\ii'}\i \ps_\ii:K_{>0}^l@>\si>>K_{>0}^l$.
From \cite{L94} it is known that
$u_{\ii,\ii'}$ is given by substraction free rational functions, so that we can
substitute $K_{>0}$ by any semifield $F$ and get a bijection $u_{\ii,\ii',F}:F^l@>\si>>F^l$.
Let $U^+_w(F)$ be the set of all $(\boc_\ii)_{\ii\in\ci(w)}\in\prod_{\ii\in\ci(w)}F^l$
such that $u_{\ii,\ii',F}(\boc_\ii)=\boc_{\ii'}$ for any $\ii,\ii'$ in $\ci(w)$.
Then $U^+(F)=\sqc_{w\in W}U^+_w(F)$ has a natural monoid structure and for $F=K_{>0}$
one can identify $U^+(F)=U^+_{\ge0}$ as monoids; for $F=\{1\}$ we have $U^+(F)=W$.

\head 2. ``Cell'' decomposition of $U^{+\t}_{\ge0}$\endhead
\subhead 2.1\endsubhead
Let $w\m w^*$ be the involution of $W$ such that $s_i^*=s_{i^*}$ for any $i\in I$.
Let $\II_*=\{w\in W;ww^*=1\}$ be the set of twisted involutions of $W$.
For any $w\in W$, $\t$ restricts to a bijection

(a) $U^+_{\ge0,w}@>\si>>U^+_{\ge0,(w^*)\i}$.
\nl
To see this we use that for $\ii=(i_1,i_2,\do,i_k)\in\ci(w)$ and
$(a_1,a_2,\do,a_k)\in K_{>0}^k$ we have
$$\t(x_{i_1}(a_1)x_{i_2}(a_2)\do x_{i_k}(a_k))=
x_{i_k^*}(a_k)x_{i_{k-1}^*}(a_{k-1})\do x_{i_1^*}(a_1).$$
In particular, for $w\in\II_*$, $\t$ defines
a bijection $U^+_{\ge0,w}@>\si>>U^+_{\ge0,w}$
whose fixed point is denoted by $U^{+\t}_{\ge0,w}$.
We see that we have a partition

(b) $U^{+\t}_{\ge0}=\sqc_{w\in\II_*}U^{+\t}_{\ge0,w}$.
\nl
The subsets $U^{+\t}_{\ge0,w}$ of $U^{+\t}_{\ge0}$ are said to be the ``cells'' of
$U^{+\t}_{\ge0}$.

\subhead 2.2\endsubhead
Let $w\in\II_*,i\in I$ be such that setting $s=s_i\in S$
we have $s^*w=ws,|w|=|s^*w|+1$. Note that $s^*w\in\II_*$.
The map in 1.4(a) restricts to a map

(a) $U^{+\t}_{\ge0,sw}\T K_{>0}@>>>U^{+\t}_{\ge0,w}$.
\nl
We show:

(b) {\it the map (a) is a bijection.}
\nl
The fact that (a) is injective follows from 1.4(b).
We prove that (a) is surjective.
Let $u'\in U^{+\t}_{\ge0,w}$. By 1.4(b) we have $u'=x_{i^*}(a)ux_i(a)$ with 
$u\in U^+_{\ge0,sw},a\in K_{>0}$.
Since $\t(u')=u'$, we have $x_{i^*}(a)\t(u)x_i(a)=x_{i^*}(a)ux_i(a)$;
we deduce that $\t(u)=u$. Thus (a) is surjective hence a
 bijection. 

\subhead 2.3\endsubhead
Let $w\in\II_*$, $i\in I$ be such that, setting $s=s_i$
we have $s^*w\ne ws,|w|=|s^*w|+1$. It follows that $s^*ws\in\II_*,|s^*ws|=|w|-2$. We define

(a) $U^{+\t}_{\ge0,s^*ws}\T K_{>0}@>>>U^{+\t}_{\ge0,w}$
\nl
by $(u,a)\m x_{i^*}(a)ux_i(a)$. This is well defined since $s^*\bul(s^*ws)\bul s=w$
(see 1.3(a)). We show:

(b) {\it the map (a) is a bijection.}
\nl
The fact that (a) is injective follows immediately from
1.2(a) since $|w|=|s^*ws|+2$.  We show that (a) is
surjective. Let $u'\in U^{+\t}_{\ge0,w}$.
Since $u'\in U^+_{\ge0,w}$ we have $u'=x_{i^*}(a)ux_i(b)$ for some 
$u\in U^+_{\ge0,s^*ws}$ and $a,b$ in $K_{>0}$.
Since $\t(u')=u'$ we have $x_{i^*}(b)\t(u)x_i(a)=x_{i^*}(a)ux_i(b)$.
We have $\t(u)\in U^+_{\ge0,s^*ws}$ since $s^*ws\in\II_*$;
using 1.2(a) we deduce that $a=b$ and $\t(u)=u$.
Thus (a) is surjective hence a bijection.

\subhead 2.4\endsubhead
Let $w\in\II_*$.
Following \cite{RS90} let $\cj(w)$ be the set of all sequences $\ii=(i_1,i_2,\do,i_k)$ in $I$
such that 
$$w=s_{i_k^*}\bul s_{i_{k-1}^*}\bul\do\bul s_{i_1^*}\bul s_{i_1}\bul s_{i_2}\bul\do
\bul s_{i_k}$$
with $k$ minimum possible. This set is nonempty; the minimum value of $k$ is denoted by
$||w||$. It is known that
$$||w||=(1/2)(|w|+\ph(w))\in\NN$$
where $\ph:\II_*@>>>\NN$ is the unique function such that $\ph(1)=0$ and such that for any
$s\in S,w\in\II_*$ with $|s^*w|=|w|-1$ we have
$\ph(w)=\ph(s^*w)+1$ if $s^*w=ws$ and $\ph(w)=\ph(s^*ws)$ if $s^*w\ne ws$. 
(If $*=1$, $\ph(w)$ is the dimension of the $(-1)$-eigenspace of $w$
on the reflection representation of $W$.)

For $\ii\in\cj(w),k$ as above, we define
$$\k_\ii:K_{>0}^k@>>>U^{+\t}_{\ge0}$$ by
$$(a_1,a_2,\do,a_k)\m x_{i_k^*}(a_k)x_{i_{k-1}^*}(a_{k-1})\do x_{i_1^*}(a_1)
x_{i_1}(a_1)x_{i_2}(a_2)\do x_{i_k}(a_k).$$
We show:

(a) {\it $\k_\ii$ defines a bijection $K_{>0}^k@>>>U^{+\t}_{\ge0,w}$.}
\nl
We argue by induction on $|w|$. If $w=1$ we have 
$U^+_{\ge0,w}=U^{+\t}_{\ge0,w}=\{1\}$ so that (a) is obvious.
Assume now that $w\ne1$. We set $s=s_{i_k}$, $\ii'=(i_1,i_2,\do,i_{k-1})$.
We set $w'=s^*w$ if $s^*w=ws$ and $w'=s^*ws$ if $s^*w\ne ws$.
Then $w'\in\II_*$, $|w'|<|w|$, $||w'||=||w||-1$  and $\ii'\in\cj(w')$.
It is enough to show that the map
$U^{+\t}_{\ge0,w'}\T K_{>0}@>>>U^{+\t}_{\ge0,w}$
given by $(u,a)\m x_{i_k^*}(a)ux_{i_k}(a)$ is a well defined bijection.
When $s^*w=ws$ this follows from 2.2; when $s^*w\ne ws$ this follows from 2.3.
This proves (a).

\subhead 2.5\endsubhead
For $w\in\II_*$, we set $U^{-\t}_{\ge0,w}=\{u\in U^-_{\ge0,w};\t(u)=u\}$.
We have a partition

(a) $U^{-\t}_{\ge0}=\sqc_{w\in\II_*}U^{-\t}_{\ge0,w}$.
\nl
For $w\in\II_*$, $\ii=(i_1,i_2,\do,i_k)\in\cj(w)$ we define
$$\k^-_\ii:K_{>0}^k@>>>U^{-\t}_{\ge0}$$ by
$$(a_1,a_2,\do,a_k)\m y_{i_k^*}(a_k)y_{i_{k-1}^*}(a_{k-1})\do y_{i_1^*}(a_1)
y_{i_1}(a_1)y_{i_2}(a_2)\do y_{i_k}(a_k).$$
Then the following analogue of 2.4(a) holds (with a similar proof):

(b) {\it $\k^-_\ii$ defines a bijection $K_{>0}^k@>>>U^{-\t}_{\ge0,w}$.}

\subhead 2.6\endsubhead
Consider the action 0.1(a) of $G$ on $G$. When an element of
the ``cell'' $U^+_{\ge0,w}$ ($w\in W$) of $U^+_{\ge0}$ is applied to an
element
of the ``cell'' $U^{+\t}_{\ge0,w'}$ ($w'\in\II_*$) of $U^{+\t}_{\ge0}$, the result is an
element of the ``cell'' $U^{+\t}_{\ge0,w\bul w'\bul w^*}$ of $U^{+\t}_{\ge0}$.
(This follows by applying 1.3(a) twice.)
We see that the action 0.1(a) induces an action of the monoid $W$ (viewed as the
indexing set of the ``cells'' of $U^+_{\ge0}$) on the set $\II_*$ 
(viewed as the indexing set of the ``cells'' of $U^{+\t}_{\ge0}$).
This action is given by:

(a) $w:w'\m w\bul w'\bul w^*$
\nl
It is remarkable that the same action appears in \cite{LV21} as the limit
when the parameter is specialized to $0$ of the action \cite{LV12} of an Iwahori-Hecke algebra
on a module with basis indexed by $\II_*$,

\subhead 2.7\endsubhead
It is clear that $U^{+\t}_{\ge0}\sub U^{+\t}(K)$. In the remainder of this subsection
we assume that $K=\RR$. We show:

(a) {\it $U^{+\t}_{\ge0,w_0}$ is an open subset of $U^{+\t}(K)$.}
\nl
Let $\ii\in\cj(w_0)$. Then $\k_\ii$ is an injective continuous map from
$K_{>0}^{||w_0||}$ to $U^{+\t}(K)$ which can be shown to be homeomorphic to
$K^{||w_0||}$. Hence $\k_\ii$ is a homeomorphism on an open subset of 
$U^{+\t}(K)$ which must be $U^{+\t}_{\ge0,w_0}$. (We have used Brouwer's theorem on
invariance of domain.) This proves (a).

Similarly,

(b) {\it $U^{-\t}_{\ge0,w_0}$ is an open subset of $U^{-\t}(K)$.}

\head 3. ``Cell'' decomposition of $G^\t_{\ge0}$\endhead
\subhead 3.1\endsubhead
For any $w,w'$ in $W$, $\t$ restricts to a bijection

(a) $G_{\ge0,w,w'}@>>>G_{\ge0,(w^*)\i,(w'{}^*)\i}$.
\nl
Indeed if $u\in U^-_{\ge0,w},u'\in U^+_{\ge0,w'},t\in T_{>0}$
then $\t(utu')=\t(u')\t(t)\t(u)$ where
$\t(u')\in U^+_{\ge0,(w'{}^*)\i}$ (see 2.1),
and similarly $\t(u)\in U^-_{\ge0,(w^*)\i}$.
In particular, for $w,w'$ in $\II_*$, $\t$ defines
a bijection $G_{\ge0,w,w'}@>>>G_{\ge0,w,w'}$
whose fixed point is denoted by $G^\t_{\ge0,w,w'}$.
We see that we have a partition

(a) $G^\t_{\ge0}=\sqc_{w,w'\text{ in }\II_*}G^\t_{\ge0,w,w'}$.
\nl
The subsets $G^\t_{\ge0,w,w'}$ of $G^\t_{\ge0}$ are said to be the ``cells'' of
$G^\t_{\ge0}$.

\subhead 3.2\endsubhead
Let $w,w'$ in $\II_*$ and $i\in I$ be such that, setting
$s=s_i\in S$, we have $s^*w=ws,|w|=|s^*w|+1$. Note that $s^*w\in\II_*$.
The map in 1.6(a) restricts to a map

(a) $G^\t_{\ge0,s^*w,w'}\T K_{>0}@>>>G^\t_{\ge0,w,w'}$.
\nl
We show:

(b) {\it the map (a) is a bijection. }
\nl
The fact that (a) is injective follows from 1.6(b).
We prove that (a) is surjective.
Let $g'\in G^\t_{\ge0,w,w'}$. By 1.6(b) we have
$g'=y_{i^*}(a)gy_i(a)$ with $g\in G_{\ge0,s^*w},a\in K_{>0}$.
Since $\t(g')=g'$ we have $y_{i^*}(a)\t(g)y_i(a)=y_{i^*}(a)gy_i(a)$;
we deduce that $\t(g)=g$. Thus (a) is surjective hence it is a
bijection. 

\subhead 3.3\endsubhead
Let $w,w'$ in $\II_*$ and $i\in I$ be such that setting
$s=s_i\in S$ we have $s^*w\ne ws,|w|=|s^*w|+1$. Note that
$s^*ws\in\II_*$, $|s^*ws|=|w|-2$. We define

(a) $G^\t_{\ge0,s^*ws,w'}\T K_{>0}@>>>G^\t_{\ge0,w,w'}$
\nl
by $(g,a)\m y_{i^*}(a)gy_i(a)$. This is well defined since $s^*\bul(s^*ws)\bul s=w$
(see 1.3(c)). We show:

(b) {\it the map (a) is a bijection. }
\nl
To prove injectivity of (a) we consider $g,g'$ in $G^\t_{\ge0,s^*ws,w'}$ and $a,a'$
in $K_{>0}$ such that $y_{i^*}(a)gy_i(a)=y_{i^*}(a')g'y_i(a')$. By 1.8(a) we have $a=a'$
hence $g=g'$. Thus (a) is injective.

We show that (a) is surjective. Let $g'\in G^\t_{\ge0,w,w'}$.
Since $g\in G_{\ge0,w,w'}$ we have $g=y_{i^*}(a)g'y_i(b)$
for some $g\in G_{\ge0,s^*ws,w'}$ and $a,b$ in $K_{>0}$. Since
$\t(g')=g'$ we have $y_{i^*}(b)\t(g)y_i(a)=y_{i^*}(a)gy_i(b)$.
We have $\t(g)\in G_{\ge0,s^*ws,w'}$ since $s^*ws\in\II_*$.
Using again 1.8(a) we see that the last equality implies $a=b$. It follows that
$\t(g)=g$. Thus (a) is a bijection.

\subhead 3.4\endsubhead
Let $w,w'$ in $\II_*$ and $i\in I$ be such that, setting
$s=s_i\in S$, we have $s^*w'=w's,|w'|=|s^*w'|+1$. Note that $s^*w'\in\II$.
The map in 1.7(a) restricts to a map

(a) $G^\t_{\ge0,w,s^*w'}\T K_{>0}@>>>G^\t_{\ge0,w,w'}$.
\nl
A proof entirely similar to that in 3.2 shows:

(b) {\it the map (a) is a bijection. }

\subhead 3.5\endsubhead
Let $w,w'$ in $\II_*$ and $i\in I$ be such that, setting
$s=s_i\in S$, we have $s^*w'\ne w's,|w'|=|s^*w'|+1$. Note that $s^*w's\in\II_*$.
We define

(a) $G^\t_{\ge0,w,s^*w's}\T K_{>0}@>>>G^\t_{\ge0,w,w'}$
\nl
by $(g,a)\m x_{i^*}(a)gx_i(a)$.
This is well defined since $s^*\bul(s^*w's)\bul s=w';$ (see 1.3(c)).
A proof entirely similar to that in 3.3 shows:

(b) {\it the map (a) is a bijection. }

\subhead 3.6\endsubhead
We consider two copies $I^1=\{i^1;i\in I\}$, $I^{-1}=\{i^{-1};i\in I\}$ of $I$.
For $w,w'$ in $\II_*$ with $||w||=k,||w'||=k'$ let $\cj(w,w')$
be the set of all sequences $\jj=(i_1^{\e_1},i_2^{\e_2},\do,i_{k+k'}^{\e_{k+k'}})$
in $I^1\sqc I^{-1}$ (here $\e_1,\do,\e_{k+k'}$ are $1$ or $-1$)
such that the subsequence consisting of terms with exponent $-1$ is obtained from
a sequence in $\cj(w)$ (by attaching the exponent $-1$)
and the subsequence consisting of terms with exponent $1$ is obtained from
a sequence in $\cj(w')$  (by attaching the exponent $1$).

Let $T_{>0}^\t=\{t\in T_{>0};\t(t)=t\}\sub G^\t_{\ge0}$.

For $\jj\in\cj(w,w')$ as above we define
$\k_\jj:K_{>0}^{k+k'}\T T_{>0}^\t@>>>G^\t_{\ge0}$ by
$$\align&(a_1,a_2,\do,a_{k+k'},t)\m \\&
x^{\e_{k+k'}}_{i_{k+k'}^*}(a_{k+k'})
x_{i_{k+k'-1}^*}^{\e_{k+k'-1}}(a_{k+k'-1})\do
x_{i_1^*}^{\e_1}(a_1)tx_{i_1}^{\e_1}(a_1)x_{i_2}^{\e_2}(a_2)\do
x_{i_{k+k'}}^{\e_{k+k'}}(a_{k+k'})\endalign$$
where we set $x_i^\e(a)=x_i(a)$ if $\e=1$, $x_i^\e(a)=y_i(a)$ if $\e=-1$. We show:

(a) {\it $\k_\jj$ defines a bijection $K_{>0}^{k+k'}\T T_{>0}^\t@>\si>>G^\t_{\ge0,w,w'}$.}
\nl
We argue by induction on $|w|+|w'|$. If $w=w'=1$ we have 
$G^\t_{\ge0,w,w'}=T_{>0}^\t$ so that (a) is obvious.

We now assume that $(w,w')\ne(1,1)$ so that $k+k'\ge1$. We set $s=s_{i_{k+k'}}\in S,
\e=\e_{k+k'}$. We set $\jj'=(i_1^{\e_1},i_2^{\e_2},\do,i_{k+k'-1}^{\e_{k+k'-1}})$.

Assume first that $\e=1$. Then $w'\ne1$. We set $w'_1=s^*w'$ if $s^*w'=w's$ and
$w'_1=s^*w's$ if $s^*w'\ne w's$. Then $w'_1\in\II_*$, $|w'_1|<|w'|$, $||w'_1||=||w'||-1$ and
$\jj'\in\cj(w,w'_1)$.
It is enough to show that the map
$G^\t_{\ge0,w,w'_1}\T K_{>0}@>>>G^\t_{\ge0,w,w'}$
given by $(g,a)\m x_{i_{k+k'}^*}(a)gx_{i_{k+k'}}(a)$ is a well defined bijection.
When $s^*w'=ws'$ this follows from 3.4; when $s^*w'\ne ws'$ this follows from 3.5.

Next we assume that $\e=-1$. Then $w\ne1$. We set $w_1=s^*w$ if $s^*w=ws$ and $w_1=s^*ws$ if
$s^*w\ne ws$. Then $w_1\in\II_*$, $|w_1|<|w|$, $||w_1||=||w||-1$ and $\jj'\in\cj(w_1,w')$.
It is enough to show that the map
$G^\t_{\ge0,w_1,w'}\T K_{>0}@>>>G^\t_{\ge0,w,w'}$
given by $(g,a)\m y_{i_{k+k'}^*}(a)gy_{i_{k+k'}}(a)$ is a well defined bijection.
When $s^*w=ws$ this follows from 3.2; when $s^*w\ne ws$ this follows from 3.3.

This completes the proof of (a).

\subhead 3.7\endsubhead
We show:

(a) $G^\t_{\ge0}\sub (G^\t)^0(K)$
\nl
(notation of 0.1). By 3.1(a) it is enough to show that 
$G^\t_{\ge0,w,w'}\sub (G^\t)^0$ for any $w,w'$ in $\II_*$. (The inclusion
$G^\t_{\ge0}\sub G(K)$ is obvious.)

Using 3.6(a) we see that it is enough to show that $T_{>0}^\t\sub (G^\t)^0$.
It is enough to show that any element $t_1\in T_{>0}^\t$ is contained in the image of
$c:T_{>0}@>>>T_{>0}$, $t\m t\t(t)$ (the converse is obvious).
We shall use the following fact:

(b) the map $T_{>0}@>>>T_{>0}$, $e\m e^2$ is an isomorphism.
\nl
Assuming (b) we write $t_1=t_2^2$ where $t_2\in T_{>0}$.
We have $(t_2\t(t_2))^2=t_2^2\t(t_2^2)=t_1\t(t_1)=t_1t_1=t_1^2$. Using the injectivity
of the map in (b) we deduce $t_2\t(t_2)=t_1$ completing the proof of (a).

To prove (b) it is enough to show that

(c) the homomorphism $K_{>0}@>>>K_{>0}$, $a\m a^2$ is an isomorphism. 
\nl
Let $a\in K_{>0}$. From the definition we have $a=b^2=(-b)^2$ for some $b\in K-\{0\}$.
Since either $b>0$ or $-b>0$, the map (c) is surjective. If $a\in K_{>0}$ satisfies
$a^2=1$ then $a=1$ or $a=-1$. But $a$ cannot be $-1$ since $-1$ is not a square in $K$.
Thus $a=1$ and the map (c) is injective hence an isomorphism.

\mpb

In the remainder of this subsection (and the next one) we assume that $K=\RR$. We have

(d) $\dim(G/H)=2||w_0||+\dim (T^\t)$
\nl
where $\dim$ is complex dimension and $T^\t=\{t\in T;\t(t)=t\}$.
An equivalent formula is $\dim(H)=|w_0|-\ph(w_0)+\dim T-\dim(T^\t)$
where $\ph(w_0)$ is as in 2.4. This follows from the definitions.
(If $G$ is almost simple, $w_0$ is central in $W$ and $*=1$ then the last equality
becomes $\dim(H)=|w_0|$.)

We show:

(e) $G^\t_{\ge0,w_0,w_0}$ is an open subset of $(G^\t)^0(K)$.
\nl
Let $\jj\in\cj(w_0,w_0)$. Then $\k_\jj$ is an injective continuous map from
$K_{>0}^{2||w_0||}\T T_{>0}^\t$ to $(G^\t)^0(K)$ which is a manifold of
real dimension $2||w_0||+\dim (T_{>0}^\t)$ (see (d))
Hence $\k_\jj$ is a homeomorphism on an open subset of $(G^\t)^0(K)$ 
which must be $G^\t_{\ge0,w_0,w_0}$. (We have used Brouwer's theorem on
invariance of domain.) This proves (e).

\subhead 3.8\endsubhead
By the exponential map $T_{>0},\t$ become an $\RR$-vector spacw and a linear involution
on it. It follows that $T_{>0}^\t$ is a cell.
Hence from 3.6 we see that:

(a) {\it for any $w,w'$ in $\II_*$, $G^\t_{\ge0,w,w'}$ is a cell of dimension $||w||+||w'||+
\dim(T_{>0}^\t)$. } and that

(b) {\it the partition 3.1(a) of $G^\t_{\ge0}$ is a cell decomposition.}

\subhead 3.9\endsubhead
As in 2.6 we see that the action 0.1(a) induces an action of the monoid $W\T W$
(viewed as the indexing set of the ``cells'' of $G_{\ge0}$) on the set $\II_*\T\II_*$ 
(viewed as the indexing set of the ``cells'' of $G^\t_{\ge0}$).
This action is given by:

(a) $(w_1,w_2):(w,w')\m (w_1\bul w\bul w_1^*,w_2\bul w'\bul w_2^*)$.

\head 4. Transition maps\endhead
\subhead 4.1\endsubhead
Now let $w\in\II_*$ with $||w||=k$. For any $\ii,\ii'$ in $\cj(w)$ we set
$v_{\ii,\ii'}=\k_{\ii'}\i \k_\ii :K_{>0}^k@>\si>>K_{>0}^k$
(a transition map). 

\subhead 4.2\endsubhead
Hu and Zhang \cite{HZ16}, \cite{HZ17} have shown (at least when $W$ is of classical type and
$*=1$) that $\cj(w)$ can be viewed as the set of vertices of a (connected) graph in which
$\ii,\ii'$ are joined when they are related by certain elementary moves which only involve a
small number of indices; these elementary moves include the standard braid moves but there are
also a small number of non-standard moves. The results of \cite{HZ16},\cite{HZ17}
were extended to $W$ of type $F_4$ in \cite{HZW} and to a general $(W,*)$ in \cite{M17}
(with the use of a computer) and in \cite{HH19} (without the use of a computer).

Note that, by the connectedness of the graph $\cj(w)$, any transition map $v_{\ii,\ii'}$ in 4.1
is a composition of transition maps corresponding to pairs $\ii,\ii'$ which form an edge of the
graph. When the edge corresponds to a standard braid move, the corresponding transition map is
of the type appearing in \cite{L94}; it involves only (substraction free) rational functions.
When the edge corresponds to a nonstandard braid move, the corresponding transition map is
of a type appearing in the next subsection.

\subhead 4.3\endsubhead
Here are examples of transition maps associated to non-standard braid moves.
(This is actually the complete list of examples associated to ``irreducible''
non-standard braid moves.) In each of the examples below we have $w=w_0$.

(i) Assume first that $W$ is of type $A_2$, $*=1$ and $I=\{1,2\}$.
There is a unique pair of inverse bijections $K^2_{>0}\lra K^2_{>0}$,
$(a,b)\lra(a',b')$ such that
$$x_1(a)x_2(b)x_2(b)x_1(a)=x_2(a')x_1(b')x_1(b')x_2(a').$$
This is a pair of inverse transition maps associated to a non-standard move
from $\ii=(1,2)$ to $\ii'=(2,1)$ in Hu-Zhang \cite{HZ16}.
Note that we must have $a'=b,b'=a$.

(ii) Assume next that $W$ is of type $A_2$, $*\ne1$ and $I=\{1,2\}$.
There is a unique pair of inverse bijections $K^2_{>0}\lra K^2_{>0}$,
$(a,b)\lra(a',b')$ such that
$$x_1(a)x_2(b)x_1(b)x_2(a)=x_2(a')x_1(b')x_2(b')x_1(a').$$
This is a pair of transition maps associated to a non-standard move
from $\ii=(1,2)$ to $\ii'=(2,1)$ in Marberg \cite{M17}.

(iii) Assume next that $W$ is of type $B_2$, $*=1$ and $I=\{1,2\}$.
There is a unique pair of inverse bijections $K^3_{>0}\lra K^3_{>0}$,
$(a,b,c)\lra(a',b',c')$ such that
$$x_1(a)x_2(b)x_1(c)x_1(c)x_2(b)x_1(a)=x_2(a')x_1(b')x_2(c')x_2(c')x_1(b')x_2(a').$$
This is a pair of transition maps associated to a non-standard move
from $\ii=(1,2,1)$ to $\ii'=(2,1,2)$ in Hu-Zhang \cite{HZ17}.

(iv) Assume next that $W$ is of type $A_3$, $*\ne1$ and $I=\{1,2,3\}$ with $1^*=3,2^*=2,3^*=1$.
There is a unique pair of inverse bijections $K^4_{>0}\lra K^4_{>0}$,
$$(a_1,a_2,a_3,a'_2)\lra(b_1,b_2,b_3,b'_2)$$
such that
$$\align&x_2(a_2)x_3(a_3)x_1(a_1)x_2(a'_2)x_2(a'_2)x_3(a_1)x_1(a_3)x_2(a_2)\\&
=x_3(b_3)x_2(b_2)x_1(b_1)x_2(b'_2)x_2(b'_2)x_3(b_1)x_2(b_2)x_1(b_3).\endalign$$
This is a pair of transition maps associated to a non-standard move
from $\ii=(2,3,1,2)$ to $\ii'=(2,3,2,1)$ in Marberg \cite{M17}).  

(v) Assume next that $W$ is of type $B_3$, $*=1$ and $I=\{1,2,3\}$ with $s_1s_2$ of order $4$ and
$s_1s_3=s_3s_1$.
There is a unique pair of inverse bijections $K^6_{>0}\lra K^6_{>0}$,
$$(a,b,c,d,e,f)\lra(a',b',c',d',e',f')$$
such that
$$\align&
   x_2(a)x_1(b)x_2(c)x_3(d)x_2(e)x_1(f)x_1(f)x_2(e)x_3(d)x_2(c)x_1(b)x_2(a)\\&
   x_1(a')x_2(b')x_1(c')x_3(d')x_2(e')x_1(f')x_1(f')x_2(e')x_3(d')x_1(c')x_2(b')x_1(a').
   \endalign$$
This is a pair of transition maps associated to a non-standard move 
from $\ii=(1,2,3,2,1,2)$ to $\ii'=(1,2,3,1,2,1)$ in Hu-Zhang \cite{HZ17}.

(vi) Assume next that $W$ is of type $D_4$, $*=1$ and $I=\{0,1,2,3\}$ with $s_1,s_2,s_3$
commuting with each other but not with $s_0$.
There is a unique pair of inverse bijections $K^8_{>0}\lra K^8_{>0}$,
$$(a_1,a_2,a_3,a_0,b_1,b_2,b_3,b_0)
\lra(a'_0,a'_1,a'_2,a'_3,b'_0,b'_1,b'_2,b'_3)$$
such that

$$\align&
x_3(a_3)x_2(a_2)x_1(a_1)x_0(a_0)x_2(b_2)x_1(b_1)x_3(b_3)
x_0(b_0)x_0(b_0)x_3(b_3)x_1(b_1)x_2(b_2)x_0(a_0)\\&
x_1(a_1)x_2(a_2)x_3(a_3)
=\\&x_0(a'_0)x_3(a'_3)x_1(a'_1)x_2(a'_2)x_0(b'_0)x_1(b'_1)
x_2(b'_2)x_3(b'_3)x_3(b'_3)x_2(b'_2)x_1(b'_1)\\&
x_0(b'_0)x_2(a'_2)x_1(a'_1)x_3(a'_3)x_0(a'_0).\endalign$$
This is a pair of inverse transition maps associated to a non-standard move 
from

$\ii=(0,3,1,2,0,1,2,3)$ to $\ii'=(3,2,1,0,2,1,3,0)$
\nl
in Hu-Zhang \cite{HZ17}.

(vii) Assume next that $W$ is of type $G_2$, $*=1$ and $I=\{1,2\}$.
There is a unique pair of inverse bijections $K^5_{>0}\lra K^5_{>0}$,
$$(a,b,c,d,e)\lra(a',b',c',d',e')$$
such that
$$\align&x_1(a)x_2(b)x_1(c)x_2(d)x_1(e)x_1(e)x_2(d)x_1(c)x_2(b)x_1(a)\\&
=x_2(a')x_1(b')x_2(c')x_1(d')x_2(e')x_2(e')x_1(d') x_2(c')x_1(b')x_2(a').\endalign$$
This is a pair of inverse transition maps associated to a non-standard move from
$\ii=(1,2,1,2,1)$ to $\ii'=(2,1,2,1,2)$ in Marberg \cite{M17}.

In the following subsections we will describe the bijections $K^n_{>0}\lra K^n_{>0}$ in
(ii)-(vi) more explicitly (in case (i) the bijection is already explicit).
The deduction of 4.7(a),4.8(a),4.9(a) from the corresponding equalities in (v),(vi) was done by
rewriting those equalities as products of matrices in a standard representation of $G$ and then using
a computer to multiply those matrices. (I thank Gongqin Li for help with programming.)

\subhead 4.4\endsubhead
In the setup of 4.3(ii),
we see by calculation that  $a+b=a'+b'$,
$a^2+2ab=b'{}^2$.
Hence
$$a'=b^2/(a+b+\sqrt{\d}),b'=\sqrt{\d}$$
where $\d=a^2+2ab>0$ and
$$a=b'{}^2/(a'+b'+\sqrt{\d'}), b=\sqrt{\d'}$$
where $\d'=a'{}^2+2a'b'>0$.

\subhead 4.5\endsubhead
In the setup of 4.3(iii), assuming that the value of the
root corresponding to $2$ on the coroot corresponding to $1$ is $-2$,
we see by calculation that 
$$b'=a+c,a'=bc^2/(a+c)^2,c'=ab(a+2c)/(a+c)^2,$$
so that $a',b',c'$ can expressed in terms of $a,b,c$ using only rational functions.
On the other hand $a,b,c$ can be expressed in terms of $a',b',c'$ as follows:
$$a= b'{}^2c'/(b'(a'+c')+\sqrt{\d}),$$
$$b=a'+c',$$
$$c=\sqrt{\d}/(a'+c'),$$
where $\d=a'b'{}^2(a'+c')>0$.

\subhead 4.6\endsubhead
In the setup of 4.3(iv),  we see by calculation that 
$$a_1+a_3=b_1+b_3,a_2+a'_2=b_2+b'_2,a_1^2a'_2=b_1^2b'_2,
(a_1+a_3)a_2+2a_1a'_2=b_1(b_2+2b'_2).\tag a$$
We try to express $(a_1,a_2,a_3,a'_2)$ in terms of $(b_1,b_2,b_3,b'_2)$.
Substituting $a_1=\a$, $a_3=b_1+b_3-\a$, $a'_2=b_1^2b'_2\a^{-2},a_2=b_2+b'_2-b_1^2b'_2\a^{-2}$ in the 
last equation in (a) we obtain
$$(b_1+b_3)b_1^2b'_2\a^{-2}-2b_1^2b'_2\a\i+b_1b'_2-b_3(b_2+b'_2)=0$$
so that
$$\a\i=(b_1+\e\sqrt{\d})/b_1(b_1+b_3) \text{ and }\a=b_1(b_1+b_3)/(b_1+\e\sqrt{\d})$$
where $\d=b_3^2+b_2b'_2{}\i b_3(b_1+b_3)>0$ and $\e=\pm1$. We have
$$a_3=b_1+b_3-b_1(b_1+b_3)/(b_1+\e\sqrt{\d})=(b_1+b_3)\e\sqrt{\d}/(b_1+\e\sqrt{\d}).$$
Since $a_3>0$, it follows that $\e=1$. We have
$$a'_2=b_1^2b'_2(b_1+\sqrt{\d})^2/b_1^2(b_1+b_3)^2=b'_2(b_1+\sqrt{\d})^2/(b_1+b_3)^2.$$
We have
$$\align&a_2=b_2+b'_2-b'_2(b_1+\sqrt{\d})^2/(b_1+b_3)^2\\&
=((b_2+b'_2)(b_1+b_3)^2-b'_2(b_1^2+2b_1\sqrt{\d}+\d))/(b_1+b_3)^2\\&
=(b_2b_1^2+b_2b_1b_3+2b'_2b_1b_3-2b_1b'_2\sqrt{\d})/(b_1+b_3)^2\\&
=\fra{(b_2b_1^2+b_2b_1b_3+2b'_2b_1b_3)^2-4b_1^2b'_2{}^2\d}
{(b_2b_1^2+b_2b_1b_3+2b'_2b_1b_3+2b_1b'_2\sqrt{\d})(b_1+b_3)^2}\\&
=\fra{b_2^2b_1^4+b_2^2b_1^2b_3^2+2b_2^2b_1^3b_3}
{(b_2b_1^2+b_2b_1b_3+2b'_2b_1b_3+2b_1b'_2\sqrt{\d})(b_1+b_3)^2}\\&
=b_2^2b_1^2/(b_2b_1^2+b_2b_1b_3+2b'_2b_1b_3+2b_1b'_2\sqrt{\d})\\&
=b_2^2b_1/(b_2b_1+b_2b_3+2b'_2b_3+2b'_2\sqrt{\d}).\endalign$$

We now try to express $(b_1,b_2,b_3,b'_2)$ in terms of $(a_1,a_2,a_3,a'_2)$. 
Substituting $b_1=\b,b_3=a_1+a_3-\b,b'_2=a_1^2a'_2\b^{-2}$ in the last equation in (a) we obtain
$$(a_2+a'_2)\b^2-((a_1+a_3)a_2+2a_1a'_2)\b+a_1^2a'_2=0$$
that is
$$\b=((a_1+a_3)a_2+2a_1a'_2+\e'\sqrt{\d'})/2(a_2+a'_2)$$
where $\d'=a_2^2(a_1+a_3)^2+4a_1a_2a'_2a_3>0$, $\e'=\pm1$. We have
$$\align&b_3=a_1+a_3-b_1=a_1+a_3-((a_1+a_3)a_2+2a_1a'_2+\e'\sqrt{\d'})/2(a_2+a'_2)\\&=
(2(a_1+a_3)(a_2+a'_2)-(a_1+a_3)a_2-2a_1a'_2+\e'\sqrt{\d'})/2(a_2+a'_2)\\&=
(a_1a_2+a_2a_3+2a'_2a_3+\e'\sqrt{\d'})/2(a_2+a'_2)\\&
=\fra{(a_1a_2+a_2a_3+2a'_2a_3)^2-\d'}{2(a_2+a'_2)(a_1a_2+a_2a_3+2a'_2a_3-\e'\sqrt{\d'})}
\\&=\fra{4a'_2{}^2a_3^2+4a_2a'_2a_3^2}{2(a_2+a'_2)(a_1a_2+a_2a_3+2a'_2a_3-\e'\sqrt{\d'})}
\\&=2a'_2a_3^2/(a_1a_2+a_2a_3+2a'_2a_3-\e'\sqrt{\d'})\\&
=2a'_2a_3^2(a_1a_2+a_2a_3+2a'_2a_3+\e'\sqrt{\d'})/ 
((a_1a_2+a_2a_3+2a'_2a_3)^2-\d')\\&
=2a'_2a_3^2(a_1a_2+a_2a_3+2a'_2a_3+\e'\sqrt{\d'})/ 
(4a'_2{}^2a_3^2  +4a_2a'_2a_3^2)\\&
=(a_1a_2+a_2a_3+2a'_2a_3+\e'\sqrt{\d'})/2(a_2+a'_2).\endalign$$
We have
$$b'_2=4a_1^2a'_2 (a_2+a'_2)^2((a_1+a_3)a_2+2a_1a'_2+\e'\sqrt{\d'})^{-2}.$$
We have
$$\align&b_2=a_2+a'_2-b'_2=a_2+a'_2-4a_1^2a'_2 (a_2+a'_2)^2((a_1+a_3)a_2+2a_1a'_2+\e'\sqrt{\d'})^{-2}\\&=
(a_2+a'_2)
\fra{((a_1+a_3)a_2+2a_1a'_2+\e'\sqrt{\d'})^2-4a_1^2a'_2 (a_2+a'_2)}
{((a_1+a_3)a_2+2a_1a'_2+\e'\sqrt{\d'})^2}\\&=
(a_2+a'_2)\T\\&
\fra{a_1^2a_2^2+a_3^2a_2^2+2a_1a_3a_2^2+4a_1a'_2a_2a_3
+2\e'((a_1+a_3)a_2+2a_1a'_2)\sqrt{\d'}+\d'}
{((a_1+a_3)a_2+2a_1a'_2+\e'\sqrt{\d'})^2}.\endalign$$
When $\e'=1$ the values of $b_1,b_2,b'_2,b_3$ (expressed in terms of $a_1,a_2,a'_2,a_3$) are $>0$.
It follows that $\e'$ must be equal to $1$.

\subhead 4.7\endsubhead
Assume that we are in the setup of 4.3(v) with the underlying reductive group of type $C_3$. We see by calculation
that 
$$\align&d=d',a+c+e=b'+e',b+f=a'+c'+f', \\&e^2f=e'{}^2f',bc+be+ef=a'b'+a'e'+c'e'+e'f',\\&
abc+abe+aef+cef+e^2f=e'(b'e'+b'f'+e'f').\tag a\endalign$$

We try to express $a',b',c',d',e',f'$ in terms of $a,b,c,d,e,f$. Setting $\a=e'$ we have

$b'=a+c+e-\a$, $a'b'=bc+be+ef-(b+f)\a$.
\nl
Setting $A=abc+abe+aef+cef$, we have $\a b'(\a+f')=A$ hence $\a(a+c+e-\a)(\a+e^2f\a^{-2})=A$
that is $(a+c+e-\a)(\a^3+e^2f)=A\a$. Thus $e'=\a$ is a root of a polynomial of degree $4$ with
coefficients rational functions in $a,b,c,d,e,f$.
Then $a',b',f',c',d'$ are rational functions in $a,b,c,d,e,f,e'$.

We now
try to express $a,b,c,d,e,f$ in terms of $a',b',c',d',e',f'$. Let $u'=a'+c'+f'$, $z=e'{}^2f'$, $w'=b'+e'$.

Setting $\b=e$ we have $f=z\b^{-2}$, $b+f=u'$ hence $b=u'-z\b^{-2}$,          
$bc+\b u'=A'$ where $A'=a'b'+e'u'$, hence $c=(A'-\b u'))/(u'-z\b^{-2})$. We have
$$\align& a=w'-\b-(A'-\b u')/(u'-z\b^{-2})\\&
=((w'-\b)(u'-z\b^{-2})-(A'-\b u'))/(u'-z\b^{-2})\\&
=(-zw'\b^{-2}+z\b\i+A'')/(u'-z\b^{-2})\endalign$$
where $A''=w'u'-A'$.

Substituting these values of $a,b,c,e,f$ in the last equation (a)
we see that $e=\b$ is a root of a polynomial of degree $2$ in
$a',b',c',e',f'$. Then $a,b,c,d,f$ are rational functions in 
$a',b',c',d',e',f',e$. 

\subhead 4.8\endsubhead
Assume that we are in the setup of 4.3(v) with the underlying reductive
group of type $B_3$. We see by calculation that 
$$\align&d=d',a+c+e=b'+e',b+f=a'+c'+f', ef=e'f',\\&
bc+be+ef=a'b'+a'e'+c'e'+e'f',\\&
2abef+2aef^2+2cef^2+ab^2c+ab^2e=e'b'((c'+f')^2+f'{}^2).\tag a\endalign$$

We try to express $a',b',c',d',e',f'$ in terms of $a,b,c,d,e,f$. Setting $\a=e'$ we have

$b'=a+c+e-\a$, $a'b'=-(b+f)\a+bc+be+ef$.
\nl
hence
$$\align&c'+f'=b+f-a'=b+f-(-(b+f)\a+bc+be+ef)/(a+c+e-\a))\\&=
((b+f)(a+c+e-\a)-(-(b+f)\a+bc+be+ef))/(a+c+e-\a)\\&=
((b+f)(a+c+e)-(bc+be+ef))/(a+c+e-\a)\\&=(ab+af+cf)/(a+c+e-\a).\endalign$$

Setting $A=2abef+2aef^2+2cef^2+ab^2c+ab^2e$ we have
$$A=\a(a+c+e-\a)((ab+af+cf)^2/(a+c+e-\a)^2+e^2f^2\a^{-2})$$
that is
$$A(a+c+e-\a)\a=\a^2(ab+af+cf)^2+(a+c+e-\a)^2e^2f^2.$$
Thus $e'=\a$ is a root of a polynomial of degree $2$ with
coefficients rational functions in $a,b,c,d,e,f$.
Then $a',b',f',c',d'$ are rational functions in $a,b,c,d,e,f,e'$.

We try to express $a,b,c,d,e,f$ in terms of $a',b',c',d',e',f'$.
Let $u'=a'+c'+f'$, $z=e'f'$, $w'=b'+e'$. Setting $\b=e$ we have $f=z\b\i$,
$b+f=u'$ hence $b=u'-z\b\i$, $bc+\b u'=A'$ where $A'=a'b'+e'u'$,
hence $c=(A'-\b u'))/(u'-z\b\i)$.
We have
$$\align& a=w'-\b-(A'-\b u')/(u'-z\b\i)\\&
=((w'-\b)(u'-z\b\i)-(A'-\b u'))/(u'-z\b\i)\\&
=(-zw'\b\i+A'')/(u'-z\b\i)\endalign$$
where $A''=w'u'-A'+z$.
Substituting these values of $a,b,c,e,f$ in the last equation (a)
we see that $e=\b$ is a root of a polynomial of degree $2$ in
$a',b',c',e',f'$. Then $a,b,c,d,f$ are rational functions in 
$a',b',c',d',e',f',e$. 

\subhead 4.9\endsubhead
Assume that we are in the setup of 4.3(vi). By calculation we see that

$$\align& a_i+b_i=a'_i+b'_i \text{ for } i=0,1,2,3;\\&
(a_ia_j+a_ib_j+b_ia_j)(a_0+b_0)+b_ib_jb_0=(a'_ia'_j+a'_ib'_j+b'_ia'_j)b'_0
\text{ for $i\ne j$ in } \{1,2,3\};\\&
2a_3b_1b_2a_0(a_0+b_0)+b_1b_2b_3a_0b_0\\&
=(a'_1a'_2a'_3+a'_1a'_3b'_2+a'_1a'_2b'_3+a'_2a'_3b'_1)a'_0b'_0
+2a'_3b'_1b'_2b'_0(a'_0+b'_0).\tag a\endalign$$
We set $\b=b'_0, \a=a_0$, $m_i=a_i+b_i=a'_i+b'_i=m'_i$.
From the first two equalities in (a) we deduce for $i\ne j$ in $\{1,2,3\}$:

(b) $b'_ib'_j=d_{ij}\b\i+m_im_j$  where $d_{ij}=b_ib_ja_0-m_im_jm_0$,

(c) $b_ib_j=d'_{ij}\a\i$ where $d'_{ij}=b'_ib'_jb'_0+m'_im'_ja'_0$

Taking product over $(i,j)$ in $(1,2),(2,3),(1,3)$ we obtain

$$\align&b'_1{}^2b'_2{}^2b'_3{}^2=(d_{12}\b\i+m_1m_2)(d_{23}\b\i+m_2m_3)(d_{13}\b\i+m_1m_3)\\&
:=\D=\d_0\b^{-3}+\d_1\b^{-2}+\d_2\b\i+\d'_3\endalign$$
where $\d_t$ are polynomials in $a_i,b_i$;
$$b_1^2b_2^2b_3^2=d'_{12}d'_{23}d'_{13}\a^{-3}:=\D'=\d'\a^{-3}$$
where $\d'=d'_{12}d'_{23}d'_{13}$.
It follows that
$$b'_1b'_2b'_3=\sqrt{\D},b_1b_2b_3=\sqrt{\D'}$$

Combining this with (b),(c) we deduce

(d) $b'_k=\sqrt{\D}/(d_{ij}\b\i+m_im_j)$, $a'_k=m'_k-\sqrt{\D}/(d_{ij}\b\i+m_im_j)$

(e) $b_k=\sqrt{\D'}/(d'_{ij}\a\i)$, $a_k=m_k-\sqrt{\D'}/(d'_{ij}\a\i)$
\nl
where $i,j,k$ is any permutation of $1,2,3$.

In the last equality in (a) the left hand side is $C:=2a_3b_1b_2a_0m_0+b_1b_2b_3a_0b_0$;
the right hand side is
equal to
$$(m_1m_2m_3-m_1b'_2b'_3-m_2b'_1b'_3-m_3b'_1b'_2+2b'_1b'_2b'_3)\b(m_0-\b)
+2(m_3-b'_3)b'_1b'_2\b m_0;$$
we substitute in it $b'_ib'_j=d_{ij}\b\i+m_im_j$, $b'_1b'_2b'_3=\sqrt{\D}$. We obtain
$$\align&m_1m_2m_3\b(m_0-\b)-(d_{12}\b\i+m_1m_2)m_3\b(m_0-\b)\\&
-(d_{23}\b\i+m_2m_3)m_1\b(m_0-\b)-(d_{13}\b\i+m_1m_2)m_3\b(m_0-\b)\\&
+2(d_{12}\b\i+m_1m_2)m_3\b m_0-2\sqrt{\D}\b^2=C.\endalign$$
We see that $\b$ is a root of a polynomial of degree $4$ with coefficients rational
functions in $a_i,b_i$.
For $k\in\{1,2,3\}$ we have $b'_k{}^2=\D/(d_{ij}\b\i+m_im_j)^2$, where $\{i,j\}=\{1,2,3\}-\{k\}$
hence $b'_k$ is a square root from a rational function in $a_i,b_i,b'_0$.
For $k\in\{0,1,2,3\}$ we have $a'_k=m_k-b'_k$ 
hence $a'_k$ is a linear function in $a_k,b_k,b'_k$.

In the last equality in (a) the right side is denoted by $C'$ and in the
left hand side we substitute
$a_3=m'_3-\sqrt{\D'}/(d'_{12}\a\i)$
$b_1b_2=d'_{12}\a\i$, $b_1b_2b_3=\sqrt{\D'}$, $b_0=m'_0-\a$. We obtain
$$2(m'_3d'_{12}\a\i-\sqrt{\D'})\a m'_0+\sqrt{\D'}\a(m'_0-\a)=C',$$
that is
$$2m'_3m'_0d'_{12}-\sqrt{\d'\a^{-3}}(2\a m'_0+\a^2)=C',$$
$$\sqrt{\d'\a^{-3}}(2\a m'_0+\a^2)=2m'_3m'_0d'_{12}-C',$$
$$\d'\a^{-3}(2\a m'_0+\a^2)^2=(2m'_3m'_0d'_{12}-C')^2,$$
$$\a+4m'_0+4m'_0{}^2\a\i=(2m'_3m'_0d'_{12}-C')^2/\d'.$$
We see that $\a$ is a root of a polynomial of degree $2$ with coefficients rational
functions in $a'_i,b'_i$.

For $k\in\{1,2,3\}$ we have $b_k^2=\D'/(d_{ij}\a\i)^2$, where $\{i,j\}=\{1,2,3\}-\{k\}$
hence $b_k$ is a square root from a rational function in $a'_i,b'_i,\a$.
For $k\in\{0,1,2,3\}$ we have $a_k=m'_k-b_k$ 
hence $a_k$ is a linear function in $a'_k,b'_k,b_k$.

\subhead 4.10\endsubhead
By the computations above we see that each of the bijections $\z:K_{>0}^n@>>>K_{>0}^n$ in
4.3(i)-(vi) is of the form $(a_1,a_2,\do,a_n)\m(a'_1,a'_2,\do,a'_n)$
(up to possibly reordering the coordinates $a_i$) where

$a'_1$ is a root of a polynomial equation of degree $N_1$ with coefficients
rational functions in $a_1,a_2,\do,a_n$;

$a'_2$ is a root of a polynomial equation of degree $N_2\le N_1$ with coefficients
rational functions in $a_1,a_2,\do,a_n,a'_1$;

$a'_3$ is a root of a polynomial equation of degree $N_3\le N_2$ with coefficients
rational functions in $a_1,a_2,\do,a_n,a'_1,a'_2$,
\nl
etc. Moreover, either all $N_i$ for $\z$ are $\le2$ or all $N_i$ for $\z\i$ are $\le2$.
Thus $\z$ has in some sense a triangular form. We expect that the same holds in case
4.3(vii).

\head 5. Examples \endhead
\subhead 5.1\endsubhead
Assume that $G=GL_n(\bK)$ with the usual pinning and $*=1$. Now
$\s=\t$ maps the matrix $(a_{ij})_{i,j\in[1,n]}$ to
$((-1)^{i+j}a_{ij})\i_{i,j\in[1,n]}$.
Thus $G^\t$ is the set of all $(a_{ij})\in G$ such that
$\sum_j(-1)^{j+k}a_{ij}a_{jk}=\d_{ik}$ for all $i,k$
and $H$ is the set of all $(a_{ij})\in G$ such that $a_{ij}=0$ whenever $i+j$ is odd
(a subgroup of $G$ isomorphic to $GL_{n/2}(\bK)\T GL_{n/2}(\bK)$ if $n$ is even or to
$GL_{(n+1)/2}(\bK)\T GL_{(n-1)/2}(\bK)$ if $n$ is odd).

If in addition we have $n=2$ then $(G^\t)^0$ consists of all
$2\T2$ matrices in $G$ with equal diagonal entries and with determinant $1$;
$G^\t-(G^\t)^0$ consists of the two diagonal matrices with
entries $1,-1$.
In this case $G^\t_{\ge0}$ consists of all $2\T2$ matrices
with entries in $K_{\ge0}$ with equal diagonal entries and
with determinant $1$.

\subhead 5.2\endsubhead
Assume that $G$ is almost simple of rank $>1$, that $w_0$ is in the centre of $W$ and that
$*=1$ hence $\t=\s$. 
Let $\car\sub\Hom(T,\CC^*)$ be the set of roots and let $ht:\car@>>>\ZZ$ be the height function
so that for $\a\in\car$ written as $\ZZ$-linear combination of simple roots, $ht(\a)$ is the
sum of coefficients in this linear combination.
Let $\fg,\fh,\ft$ be the Lie algebras of $G,H,T$. Then $\fh$ is the subspace of $\fg$
spanned by $\ft$ and by the root subspaces corresponding to roots of even height.
This is a simple Lie algebra for which the set of simple roots with respect to
$\ft$ consists of the $\dim\ft-1$ roots of height $2$ in $\car$ and the unique root of height
$-(h-2)$ in $\car$. (Here $h$ is the Coxeter number.)
For example, if $G$ is of type $D_{2n},E_7,E_8$, then $H$ is of type $D_n\T D_n, A_7,D_8$.

\subhead 5.3\endsubhead
In this subsection we replace $G$ by $G\T G$ with the pinning
induced from that of $G$; we define $\o:G\T G@>>>G\T G$
by $(g,g')\m(g',g)$. Then $\t(g,g')=(\s(g'),\s(g))$ where
$\s$ refers to $G$. We have $(G\T G)^\t=\{(g,g')\in G\T G;g'=\s(g)\}$.
We have $(G\T G)_{\ge0}=G_{\ge0}\T G_{\ge0}$ and
$(G\T G)_{\ge0}^\t=\{(g,g')\in G_{\ge0}\T G_{\ge0};g'=\s(g)$.
This can be identified with $G_{\ge0}$ of 0.4 by
$(g,g')\m g$. In this way the theory of total positivity in \cite{L94}
becomes a special case of the theory in this paper.

\head 6. Passage to zones\endhead
\subhead 6.1\endsubhead
Let $F$ be a semifield. Now $\t$ acts naturally as an involutive antiautomorphism of the monoid $U^+(F)$
in 1.9 and one could define $U^{+\t}(F)$ as the fixed point set of $\t:U^+(F)@>>>U^+(F)$.
(One could give a similar definition for $G^\t(F)$.)
But with this definition it is not clear how to parametrize $U^{+\t}(F)$ or $G^\t(F)$
by a union of pieces of the
form $F^n$ when $F$ is other than $K_{>0}$ or $\{1\}$.

In the remainder of this section $U^{+\t}(F)$ will not refer to the above definition.
Instead we will try to redefine it for certain $F$ and certain $G$ using the method of passage to
zones in \cite{L94}.

We now assume that $K$ is the field of Puiseux series in a variable $\x$ with real coefficients
(by a theorem of Newton and Puiseux, this field is real closed). Any $x\in K-\{0\}$ is of the form
$\sum_{m\in e(x)}a_m\x^m$ where $a_m\in \RR-\{0\}$ and
$e(x)$ is a nonempty subset of $\QQ$ such that $ne(x)\sub\ZZ$ for
some $n\in\{1,2,\do\}$ and $a+e(x)\sub\QQ_{\ge0}$ for some $a\in\QQ$; note that $e(x)$ has
a well defined smallest element $v(x)\in\QQ$.
Note that $K_{>0}$ is the set of all $x\in K-\{0\}$ such that $a_{v(x)}(x)\in\RR_{>0}$.
We regard $\QQ$ as a semifield with the product of $a,b$ being $a+b$ and the sum of $a,b$ being
$\min(a,b)$. Then $v:K_{>0}@>>>\QQ$ is a semifield homomorphism. 

Let $k\in\NN$. We define a {\it zone} of $K_{>0}^k$ to be any fibre of the map
$v^k:K_{>0}^k@>>>\QQ^k$ given by $(x_1,x_2,\do,x_k)\m(v(x_1),v(x_2),\do,v(x_k))$.

\subhead 6.2\endsubhead
In this section (until the end of 6.6) we assume that $G$ is almost simple and that either

(i) $*=1$ and $G$ is of type $A$, or

(ii) $*=1$ and $G$ is of type $B_2$, or

(iii) $*\ne1$ (so that $G$ is of type $A,D$ or $E_6$).
\nl
Now let $w\in\II_*$ with $||w||=k$. Let
${}'U^{+\t}_{\ge0,w}$ be the set of all $(\boc_\ii)_{\ii\in\cj(w)}\in\prod_{\ii\in\cj(w)}K_{>0}^k$
such that $v_{\ii,\ii'}(\boc_\ii)=\boc_{\ii'}$ (notation of 4.1) for any $\ii,\ii'$ in $\cj(w)$.

We define $\k_w:{}'U^{+\t}_{\ge0,w}@>>>U^{+\t}_{\ge0,w}$ by

$\k_w((\boc_\ii)_{\ii\in\cj(w)})=\k_\ii(\boc_\ii)$ for some/any $\ii\in\cj(w)$ (notation of 2.4).
\nl
From 2.4(a) it follows that

(a) $\k_w$ is a bijection.
\nl
For $\ii,\ii'$ in $\cj(w)$, the following holds:

(b) $v_{\ii,\ii'}:K_{>0}^k@>>>K_{>0}^k$ maps any zone to a zone.
\nl
We can assume that $\ii,\ii'$ are related by
an elementary move (see 4.2). If this elementary move involves a standard braid move then (b) follows from the
results in \cite{L94}; if it involves a non-standard move then we use the formulas in  4.3(i) (in case
(i)), 4.5 (in case (ii)) or 4.4, 4.6 (in case (iii)).

From (b) we deduce that $v_{\ii,\ii'}:K_{>0}^k@>>>K_{>0}^k$ induces a map from the set of zones
in $K_{>0}^n$ to itself that is a map (necessarily a bijection) $\bar v_{\ii,\ii'}:\QQ^k@>>>\QQ^k$.

We define $U^{+\t}_w(\QQ)$ to be the set of all
$(\boc_\ii)_{\ii\in\cj(w)}\in\prod_{\ii\in\cj(w)}\QQ^k$
such that $\bar v_{\ii,\ii'}(\boc_\ii)=\boc_{\ii'}$ for any $\ii,\ii'$ in $\cj(w)$.

We define a map

$U^{+\t}_{\ge0,w}@>>>U^{+\t}_w(\QQ)$
\nl
by $u\m\bar u$ where the $\ii$-coordinate of $\bar u$ is $v^k$ applied to the $\ii$-th coordinate
of $\k_w\i(u)$. The fibres of this map are called the zones of $U^{+\t}_{\ge0,w}$. 
Taking disjoint union over $w\in\II_*$ we obtain a (surjective) map

$U^{+\t}_{\ge0}@>>>U^{+\t}(\QQ):=\sqc_{w\in\II_*}U^{+\t}_w(\QQ)$.

In \cite{L94, \S9} a definition of zones of $U^+_{\ge0}$ is given, so that the set of zones
which may be denoted $U^+(\QQ)$ is defined. (Actually in {\it loc.cit.} a subfield of $K$
on which $v$ has integer values is used instead of $K$.) From the definition we see that the
natural action of $U^+_{\ge0}$ on $U^{+\t}_{\ge0}$ induces by passage to zones an action of
$U^+(\QQ)$ on $U^{+\t}(\QQ)$.

\subhead 6.3\endsubhead
Similarly for any $(w,w')\in\II_*\T\II_*$ with $||w||=k,||w'||=k'$ we can define a partition of
$G^\t_{\ge0,w,w'}$ into zones. For this we use the various parametrizations $\k_\jj$ ($\jj\in\cj(w,w')$
of $G^\t_{\ge0,w,w'}$ in 3.6). We must show that any
$\jj,\ti\jj$ 
in $\cj(w,w')$ can be joined by a
a sequence $\jj=\jj_1,\jj_2,\do,\jj_u=\ti\jj$ in $\cj(w,w')$ so that
for any two consecutive terms $\jj_h,\jj_{h+1}$ in this sequence the composition
$\k_{\jj_{h+1}}\i\k_{\jj_h}: K_{>0}^{k+k'}\T T_{>0}^\t@>>>K_{>0}^{k+k'}\T T_{>0}^\t$
is of the type considered in \cite{L94, 1.3} or of the type
considered for $U^{+\t}_{\ge0}$ in 6.2 or of the analogous type
for $U^{=\t}_{\ge0}$.
It follows that $\k_{h+1}\i\k_h: K_{>0}^{k+k'}\T T_{>0}^\t@>>>K_{>0}^{k+k'}\T T_{>0}^\t$
maps any zone of $K_{>0}^{k+k'}$ times $T_{>0}^\t$ to a zone of $K_{>0}^{k+k'}$ times $T_{>0}^\t$.

Let $G^\t_{w,w'}(\QQ)$ be the set of zones of $G^\t_{\ge0,w,w'}$.
We set
$$G^\t(\QQ)=\sqc_{(w,w')\in\II_*\T\II_*}G^\t_{w,w'}(\QQ).$$
Now the zones of $G_{\ge0}$ are defined as in \cite{L94}. 
(Actually in {\it loc.cit.} a subfield of $K$
on which $v$ has integer values is used instead of $K$.) The set of zones of $G_{\ge0}$
is denoted by $G(\QQ)$; it inherits a monoid structure from that of $G_{\ge0}$.
The action of the monoid $G_{\ge0}$ on $G^\t_{\ge0}$ induces an action of the monoid
$G(\QQ)$ on $G^\t(\QQ)$.

\subhead 6.4\endsubhead
Let $w\in\II_*$ with $||w||=k$.
The collection of parametrizations
$\k_\ii:K_{>0}^k@>>>U^{+\t}_{\ge0,w}$ for various $\ii\in\cj(w)$
is something slightly more general than what in \cite{L19} (and going back to \cite{L94})
was called a ``positive structure''. Namely the compositions
$\k_{\ii'}\i\k_\ii:K_{>0}^k@>>>K_{>0}^k$ with $\ii,\ii'$ in $\cj(w)$
are of the form $(a_1,a_2,\do,a_k)\m(a'_1,a'_2,\do,a'_k)$
where each $a'_i$ is obtained from
$a_1,a_2,\do,a_k$ by using  a succession of the following operations:
addition, multiplication, division and extracting of a square root, not necessarily in this order.
(This last operation was not allowed in \cite{L19}.)

However in the case 6.2(i) we have just the old type of  positive structure.

\subhead 6.5\endsubhead
We will now define $U^{+\t}(F)$ for $F$ as in (i),(ii) below:

(i) $F=K_{2,>0}=K_2\cap K_{>0}$
where $K_2=\{x\in K-\{0\};2^ce(x)\sub\ZZ\text{ for some }c\in\NN\}\sqc\{0\}$;

(ii) $F=\ZZ[1/2]=\{q\in\QQ; 2^cq\in\ZZ\text{ for some }c\in\NN\}$.
\nl
Note that $K_{2,>0}$ (resp. $\ZZ[1/2]$) is a sub-semifield of $K_{>0}$ (resp. of $\QQ$).

Let $w\in\II_*$ with $||w||=k$. Let
$U^{+\t}_w(K_{2,>0})$ be the set of all $(\boc_\ii)_{\ii\in\cj(w)}\in\prod_{\ii\in\cj(w)}K_{2,>0}^k$
such that $v_{\ii,\ii'}(\boc_\ii)=\boc_{\ii'}$ (notation of 4.1) for any $\ii,\ii'$ in $\cj(w)$.
This is well defined since $v_{\ii,\ii'}$ restricts to a bijection $K_{2,>0}^k@>>>K_{2,>0}^k$
(by results in 4.4-4.6).
We define $U^{+\t}(K_{2,>0})=\sqc_{w\in\II_*}U^{+\t}_w(K_{2,>0})$;
it has a natural action of the monoid $U^+(K_{2,>0})$.

We define $U^{+\t}_w(\ZZ[1/2])$ to be the set of all
$(\boc_\ii)_{\ii\in\cj(w)}\in\prod_{\ii\in\cj(w)}(\ZZ[1/2])^k$
such that $\bar v_{\ii,\ii'}(\boc_\ii)=\boc_{\ii'}$ for any $\ii,\ii'$ in $\cj(w)$.
(Note that $\bar v_{\ii,\ii'}:\QQ^k@>>>\QQ^k$ restricts to a map
$(\ZZ[1/2])^k@>>>(\ZZ[1/2])^k$.)
We define $U^{+\t}(\ZZ[1/2])=\sqc_{w\in\II_*}U^{+\t}_w(\ZZ[1/2])$; it has a natural action of the
monoid $U^+(\ZZ[1/2])$.

We can define in a similar way $G^\t(F)$ for $F$ as in (i) or (ii).

\subhead 6.6\endsubhead
We regard $\ZZ$ as a sub-semifield of $\QQ$. It is known that the set $U^+(\ZZ)$ is
closely related to the parametrization of the canonical basis in \cite{L90}. It would be
interesting to find an analogous interpretation of $U^{+\t}(\ZZ[1/2])$.

\subhead 6.7\endsubhead
We expect that $U^{+\t}(\QQ),G^\t(\QQ)$ can be defined by the method of
6.2, 6.3 without the assumption in 6.2. This should follow from a better
understanding of the transition maps in 4.3(v)-(vii).

\enddocument